\documentclass[11pt]{amsart}
\usepackage{verbatim, amscd, epsfig,amsmath,amsthm,amstext,amssymb,
}
\theoremstyle{plain}
\newtheorem*{theorem*}{Theorem}
\newtheorem{theorem}{Theorem}[section]
\newtheorem{lemma}[theorem]{Lemma}
\newtheorem{cor}[theorem]{Corollary}

\theoremstyle{definition}

\theoremstyle{remark}
\newtheorem{rem}[theorem]{Remark}

\numberwithin{equation}{section}

\renewcommand{\Im}{{\rm im}\,}
\newcommand{\del}{\partial}
\newcommand{\delbar}{\bar{\del}}

\newcommand{\R}{\mathbb{ R}}
\newcommand{\C}{\mathbb{ C}}

\newcommand{\Zz}{\mathcal{Z}}

\renewcommand{\H}{\mathbb{ H}}

\renewcommand{\P}{\mathbb{ P}}
\newcommand{\HP}{\H\P}
\newcommand{\CP}{\C\P}

\newcommand{\trivial}[1]{\underline{\H}^{#1}}
\newcommand{\dbar}{{\bar{\partial}}}
\newcommand{\invers}{^{-1}}
\DeclareMathOperator{\End}{End}
\DeclareMathOperator{\Hom}{Hom}
\DeclareMathOperator{\pr}{pr}
\DeclareMathOperator{\ord}{ord}
\DeclareMathOperator{\Gl}{GL}
\DeclareMathOperator{\tr}{tr}

\DeclareMathOperator{\Id}{Id}

\newcommand{\delbari}{\nabla^{(0,1)}}

\newcommand{\Rr}{{\mathcal R}}
\newcommand{\K}{{\mathcal K}}

\newcommand{\plusproj}{\frac{1-IS}2}
\newcommand{\minusproj}{\frac{1+IS}2}

\begin{document}
\title
{Willmore tori in the $4$--Sphere with nontrivial normal bundle}
\author{K. Leschke, F. Pedit, U. Pinkall}

\address{ Katrin Leschke, Franz Pedit\\ 
Department of Mathematics \\
University of Massachusetts\\
Amherst, MA 01003, USA
}
\address{Ulrich Pinkall\\
Fachbereich Mathematik\\ Technische Universit\"at Berlin\\
Str. des 17.Juni 135\\
D-10623 Berlin}
\email{leschke@gang.umass.edu, franz@gang.umass.edu, pinkall@math.tu-berlin.de}

\maketitle

\section{Introduction}
The study of Willmore surfaces which are critical points for the
bending energy $\int H^2$, where $H$ is the mean curvature, goes at
least back to Blaschke's school in the $1920$'s. About $40$ years
later Willmore \cite{Willmore} reintroduced the problem and asked to
find the minimizers for the bending energy, nowadays called {\em
  Willmore energy}, over compact surfaces of fixed genus. He showed
that the round sphere is the minimum over genus zero surfaces and
formulated the conjecture that the minimum over tori is given by the
Clifford torus with Willmore energy $2\pi^2$. In the $1980$'s Bryant
\cite{Bryant} classified all Willmore spheres in 3-space as inverted
minimal spheres with planar ends in $\R^3$. Subsequently, Ejiri
\cite{Ejiri} and recently Montiel \cite{montiel} proved an analogous
result for Willmore spheres in 4-space: in addition to inverted
mininal spheres in $\R^4$ also twistor projections to $S^4$ of
rational curves in $\C\P^3$ occur.

The case of Willmore tori is more involved: there are examples
constructed by integrable system methods which are neither inverted
minimal surfaces nor twistor projections of elliptic curves
\cite{Pin_Hopf}, \cite{FerPed}, \cite{bab&bob}. By now there is a
reasonable understanding of how to construct all Willmore tori in $3$
and $4$-space from theta functions on finite genus Riemann surfaces,
the {\em spectral curves} \cite{S4}, \cite{schmidt}. In fact, the
recent preprint \cite{schmidt} by Schmidt seems to go some way
towards proving the Willmore conjecture.

An important aspect of the theory of Willmore surfaces is its
connection to the theory of harmonic maps. The {\em conformal Gau{\ss}
  map} or {\em mean curvature sphere congruence} of a Willmore surface
is a harmomic map whose energy is equal to the Willmore energy.  This
relationship between Willmore surfaces and harmonic maps becomes even
more pronounced in the description of surface theory via quaternionic
holomorphic geometry \cite{coimbra}. In this setting the theory of
Willmore surfaces in $S^4$ shows a close resemblance to the theory of
harmonic maps into $S^2$. 

A classical result of Eells and Wood \cite{eells&wood} states that a
harmonic map $f:M\to S^2$ from a compact Riemann surface $M$ of degree
$|\deg f|>\tfrac{1}{2}\deg K$, where $K$ denotes the canonical bundle
of $M$, is holomorphic or antiholomorphic. If $M$ is a torus then only
degree zero harmonic maps are nonholomorphic, and these are the
Gau{\ss} maps of constant mean curvature tori in $\R^3$. Such harmonic
maps are constructed by integrable systems methods and are given by
theta functions on hyperelliptic Riemann surfaces, the spectral curves
of the harmonic torus. Therefore, at least for tori, one can view the
result of Eells and Wood as a criterion to distinguish the trivial
holomorphic case from the more involved integrable system case.

In view of the close resemblance between harmonic maps into $S^2$ and
Willmore surfaces in $S^4$, we expect a similar criterion to hold for
Willmore surfaces: under which conditions does a Willmore surface in
$S^4$ come from a twistor projection of a holomorphic curve in
$\C\P^3$ or a minimal surface in $\R^4$, i.e., is given by
holomorphic data?
\begin{theorem*}
  Let $f:T^2\to S^4$ be a Willmore torus in $S^4$ with nontrivial
  normal bundle. Then $f$ comes from a twistor projection of an
  elliptic curve in $\C\P^3$ or from a minimal torus with planar ends in
  $\R^4$.
\end{theorem*}
In fact, we conjecture the following more general result for any
compact Willmore surface $f:M\to S^4$: if the normal bundle degree $v$
satisfies $|v|>2\deg K$, then the Willmore surface comes from a
twistor projection of a holomorphic curve in $\C\P^3$ or from a
minimal surface in $\R^4$. Of course, this is an exact analog of the
above mentioned result by Eells and Wood for harmonic maps into $S^2$.
We notice that in case $M$ has genus at most $1$ this conjecture is
true: for Willmore spheres it is the result by Ejiri and Montiel, and
for Willmore tori it is the theorem stated above. Moreover, if $f$ is
minimal in $S^4$ or, more generally, if $f$ has a dual Willmore surface then by
Theorem \ref{t:AQ=0} the conjecture holds for any genus.

The theorem by Eells and Wood follows from the fact that the
$(1,0)$--part of the derivative of a harmonic map into $S^2$ is
holomorphic together with a degree calculation. For Willmore surfaces
such a computation can also be done, but turns out to be 
insufficient for proving the theorem. The additional ingredient needed is a
detailed study of the monodromy of the associated family of Willmore
surfaces. Our model for the M\"obius geometry of $S^4$ is the
quaternionic projective line $\H\P^1$ on which the M\"obius group acts
by $\text{Gl}(2,\H)$. The associated family of Willmore surfaces is
described by an $S^1$--family of flat connections with
$\text{Gl}(2,\H)$ monodromy. Nontrivial normal bundle together with
the Pl\"ucker formula imply that, over a torus, this loop of monodromy
representations has all of its eigenvalues equal to $1$.  In case the
monodromy is trivial, the Willmore torus comes from a twistor
projection. The only other possibility is translational monodromy, in
which case the Willmore surface is an inverted minimal torus in $\R^4$
with planar ends.

In terms of spectral curves our result can also be given the following
interpretation: a Willmore torus with nontrivial normal bundle is
known to the extent one understands elliptic curves in $\CP^3$ and
minimal tori with planar ends in $\R^4$, both of which are given by
elliptic functions.  For a Willmore torus with trivial normal bundle,
which is not an inverted minimal torus in $\R^3$, the monodromy
representation of the family of flat connections has non--constant
eigenvalues. In this case, one can associate to the Willmore torus its
spectral curve, namely the Riemann surface defined by the eigenvalues
of the monodromy in dependence of the complexified loop parameter
\cite{S4}. The Willmore torus is then parameterized by theta functions
on the spectral curve, a topic which we will return to in a forthcoming
paper. 

\section{Preliminaries and degree estimates}
Before describing our setup it will be helpful to collect
some of the basic notions concerning the theory of 
quaternionic vector bundles over
Riemann surfaces \cite{Klassiker}. A quaternionic vector bundle
$W$ with complex structure $J$ over a Riemann surface $M$ decomposes into
$W = W_+ \oplus W_-$, where $W_\pm$ are the $\pm i$--eigenspaces of
$J$. By restriction $J$ induces complex structures on $W_\pm$ and $W_-
= W_+ j$ gives a complex linear isomorphism between $W_+$ and $W_-$.
The degree of the quaternionic bundle $W$ with complex structure $J$
is then defined as the degree of the underlying complex vector bundle
\begin{equation}
\label{e:degree}
\deg W := \deg W_+\,,
\end{equation}
which is half of the usual degree of $W$ viewed as a complex bundle.

Given two
quaternionic bundles $W$ and $\tilde W$ with complex structures $J$
and $\tilde J$ the complex linear homomorphisms $\Hom_+(W,\tilde W)$
are complex linearly isomorphic to $\Hom_\C(W_+,\tilde W_+)$. On the
other hand, the complex antilinear homomorphisms $\Hom_-(W,\tilde W)$
are complex linearly isomorphic to $\Hom_+(\bar{W}, \tilde W)$, where
the complex structure on a homomorphism bundle is induced by the
target complex structure. 

A {\em quaternionic holomorphic} structure on the vector bundle $W$
with complex structure $J$ is given by a quaternionic linear operator
\begin{equation}
\label{e:quat_hol_structure}
\delbar+Q:\Gamma(W)\to\Omega^{0,1}(W)=\Gamma(\bar{K}W)\,.
\end{equation}
Here $\delbar=\delbar\oplus\delbar$ is the double of a complex holomorphic structure on $W_+$ and $Q\in\Omega^{0,1}(\End_{-}(W))$ is a $(0,1)$--form with values in
complex antilinear endomorphisms of $W$. The quaternionic vector space of holomorphic sections of $W$ is denoted by
\[
H^{0}(W)=H^{0}(W,\delbar+Q)=\ker (\delbar+Q)
\]
and is finite dimensional for compact $M$. The $L^2$--norm 
\[
\mathcal{W}(W)=\mathcal{W}(W,\delbar+Q)=2\int_M <Q\wedge *Q>
\]
of $Q$ is called the {\em Willmore energy} of the holomorphic bundle $W$
where $<\,,\,>$ denotes the trace pairing on $\End(W)$. The special case $Q=0$,
for which $\mathcal{W}(W)=0$, describes (doubles of) complex holomorphic bundles $W=W_+\oplus W_+$.
A typical example of a
quaternionic holomorphic structure arises from the $(0,1)$--part $\nabla''$ of a
quaternionic connection $\nabla$ on $W$. 


Now let $f: M \to S^4$ be a conformal map of the Riemann surface $M$. We
model the M\"obius geometry of $S^4$ by the projective geometry of the
quaternionic projective line $\HP^1$. Therefore, the map $f$
corresponds to the line subbundle $L\subset V$ with $L_p = f(p)$,
where $V$ is the trivial $\H^2$--bundle over $M$. Its differential
$df$ corresponds to the $\Hom(L, V/L)$--valued 1--form
\[
\delta=\pi\nabla|_L\,,
\]
where $\pi: V \to V/L$ is the canonical projection and $\nabla$
denotes the trivial connection on $V$. A 2--sphere in $S^4$ is given
by an endomorphism $S\in\End(\H^2)$ with $S^2=-1$: points on the
2--sphere correspond to fixed lines of $S$. We denote by $\Zz$ the space
of oriented 2--spheres in $S^4$.  A sphere congruence $S:M \to\Zz$ is
thus a complex structure on $V$.

Given such a complex structure, we can decompose the tri\-vial
con\-nection in\-to $S$--com\-muting and anti\-commuting parts
\begin{equation}
\label{eq:nabla_decompose}
\nabla  = \hat\nabla+ A + Q\,,
\end{equation}
where $\hat\nabla$ is a complex connection, and $-2*A$ and $2*Q$ are
the $(1,0)$ and $(0,1)$--parts of 
\begin{equation}
\label{eq: nabla_S}
\nabla S = 2(*Q-*A)\,.
\end{equation}
By construction, $A\in\Gamma(K\End_-(V))$ and $Q\in\Gamma(\bar K
\End_-(V))$, i.e., $*A  = SA = -AS$ and $*Q=-SQ = QS$.

Among all sphere congruences the \emph{mean curvature sphere
congruence} $S: M \to\Zz$, also called the \emph{conformal Gau{\ss}
map}, of $f$ is characterized by the following properties
\cite{coimbra}:
\renewcommand{\labelenumi}{(\roman{enumi})}
\begin{enumerate}
\item \label{item:1} The sphere $S(p)$ passes through $f(p)$ for $p\in
M$, i.e., $SL=L$.
\item \label{item:2} The sphere $S(p)$ is tangent to $f$ at $p$ for
$p\in M$, i.e., $*\delta=S\delta=\delta S$.
\item \label{item:3} The sphere $S(p)$ has the same mean curvature
vector as $f$ at $p$ for $p\in M$, i.e,    $A V\subset\Omega^1(L)$,
or, equivalently, $ Q|_L=0$.
\end{enumerate}
In general, a conformal map $f: M \to S^4$ has a mean curvature sphere
congruence only along immersed points. In the sequel, we will always assume
that $f$ has a mean curvature sphere congruence which is certainly the
case when $f$ is immersed.

Note that (ii) implies that $\delta$ is a $(1,0)$--form with values in
the complex linear homomorphisms, i.e.,
$\delta\in\Gamma(K\Hom_+(L,V/L))$.  The complex connection
(\ref{eq:nabla_decompose}) decomposes into $(1,0)$ and $(0,1)$--parts
\begin{equation}
\label{eq:complex_connection}
\hat\nabla = \hat\nabla' + \hat\nabla'' =: \partial + \delbar
\end{equation}
and $\delbar$ stabilizes $L$ and therefore also $V/L$: from (ii) and
(iii) we see that $\pi\delbar|_L = \delta'' =0$. Thus, $L$ and $V/L$
are (doubles of) complex holomorphic line bundles and by (iii)
\begin{equation}
\label{eq:delta_holomorphic}
\delbar\delta = \pi\delbar\partial|_L = \pi(\partial\delbar +
R^{\hat\nabla})|_L = \delta\delbar\,,
\end{equation}
where the curvature $R^{\hat\nabla} = -(A\wedge A + Q\wedge Q)$ of
$\hat\nabla$ stabilizes $L$. This shows that $\delta$ is a holomorphic
section $\delta\in H^0(K\Hom_+(L, V/L))$ and, using \eqref{e:degree}, we obtain
\begin{equation}
\label{eq:ord_delta}
\ord\delta = \deg K + \deg V - 2\deg L\,.
\end{equation}
If $f$ is immersed, $\delta$ has no zeros and therefore
\begin{equation}
\label{eq:degV}
\deg V = 2\deg L - \deg K\,.
\end{equation}
The tangent bundle of $\HP^1$ splits into
\[
f^*(T\HP^1)= \Hom(L,V/L) =  \Hom_+(L,V/L)\oplus  \Hom_-(L,V/L)\,,
\]
where $\Hom_+(L,V/L) \supseteq \delta(TM)$ and $\Hom_-(L,V/L) =$
$\Hom_+(\bar L, V/L)$ extend the tangent bundle and the normal bundle
of $f$ across the branch points.  Therefore, the normal bundle degree
$v$ of $f$ calculates to
\begin{equation}
\label{eq:normalbundle_degree}
v = \deg V\,.
\end{equation}
Up to now, our discussion dealt with conformal maps $f: M \to
S^4$ and their mean curvature sphere congruences. In case $f$ is a
Willmore surface, we will be able to derive further degree relations.
The Willmore functional of a conformal map \cite{Klassiker} is given
by
\[
\mathcal{W}(f) = 2 \int_M <A\wedge *A>\,,
\]
which, up to topological terms, is the Willmore energy of the quaternionic holomorphic structure $\delbar +\pi Q$ on $V/L$. The Euler
Lagrange equation \cite{coimbra} of this functional is 
\begin{equation}
\label{eq:Willmore}
d^\nabla*A =0\quad \text{ or, equivalently } \quad  d^\nabla*Q=0\,,
\end{equation}
where the latter can be seen by differentiating (\ref{eq: nabla_S}).
For degree computations it is necessary to interpret $A$ and $Q$ as
complex holomorphic bundle maps. From (\ref{eq:nabla_decompose}) we
obtain
\[
d^\nabla * A = d^{\hat\nabla} * A + [A\wedge *A] + [Q\wedge *A] = 
d^{\hat\nabla}*A = Sd^{\hat\nabla} A\,,
\]
and similarly,
\[
d^\nabla *Q = -S d^{\hat\nabla} Q\,,
\]
where $[A\wedge *A]=0$ by symmetry and $[Q\wedge *A]=0$ by type
considerations. Therefore, viewing $A\in \Gamma(K\Hom_+(\bar V, V))$
and $Q\in\Gamma(K\Hom_+(V,\bar V))$, equations (\ref{eq:Willmore}) are
equivalent to
\[
\delbar A = 0 \quad  \text{ and }  \quad \delbar Q =0\,,
\]
which means that $A\in H^0(K\Hom_+(\bar V, V))$ and $Q\in
H^0(K\Hom_+(V,\bar V))$. From (iii) we see that $A$ and $Q$ have at
most rank 1, and hence there exist holomorphic subbundles $\tilde L,
\hat L \subset \bar V$, the \emph{forward} and \emph{backward
  B\"acklund transforms} \cite{coimbra} of $f$, such that
\begin{equation}
\label{eq:baecklund}
\tilde L \subseteq  \ker A \quad
 \text{ and  }\quad \hat L \supseteq  \Im Q\,.
\end{equation}
If $A\not=0$ and $Q\not=0$ the forward and backward B\"acklund
transforms are again conformal maps into $S^4$, but their mean
curvature sphere congruences may not extend into their branch points.
In case $\tilde{L}=\hat{L}$, i.e., if $AQ=0$, we will see below that the conformal map $\tilde{L}$ has mean curvature sphere congruence $-S$ and is therefore a dual Willmore surface to $f$. Now
\begin{gather}
\label{eq:A}
A\in H^0(K\Hom_+(\bar V/\tilde L, L)) \nonumber \\ 
Q\in  H^0(K\Hom_+(V/ L, \hat L)) \\
AQ \in H^0(K^2\Hom_+(V/L,L) ) \nonumber
\end{gather}
define holomorphic bundle maps between complex
holomorphic line bundles.
Therefore, the order of zeros of $A$ and $Q$  calculate to
\begin{eqnarray}
  \ord A &=& \deg K + \deg L + \deg V  - \deg \tilde L
\nonumber \\
&=&
 3 \deg L - \deg\tilde L  + \ord \delta\,, \label{eq:ord_A} 
\\  
\ord Q &=& -\deg \hat L - \deg V + \deg L + \deg K \nonumber \\
& =&
2 \deg K - \deg L  -\deg \hat  L - \ord \delta\,, \label{eq:ord_Q}  
\end{eqnarray}
where we used \eqref{eq:ord_delta}.  Moreover, if $AQ\not=0$ then
\begin{equation}
\ord AQ = 3 \deg K  - \ord \delta\,.
 \label{eq:ord_AQ}
\end{equation}
If $A = 0$ or $Q=0$ then the Willmore surface $f$ comes from
holomorphic data \cite{Klassiker}: in the former case $f$ and in the
latter case the Willmore surface $f^\perp$, given by the line bundle
$L^\perp\subset V^*$, is the twistor projection of a holomorphic curve
$g: M \to\CP^3$.

\begin{theorem}
\label{t:AQ=0}
Let $f: M\to S^4$ be a compact Willmore surface with normal bundle degree 
\[
|v| > 2\deg K
\] 
and $AQ= 0$, i.e., admitting a dual Willmore surface.  
Then, either $f$ or $f^\perp$ is a twistor projection of
a holomorphic curve in $\CP^3$, or $f$ is an inverted minimal surface
in $\R^4$.
\end{theorem}
\begin{rem} 
  If $M=S^2$ then $\deg K = -2$ and the hypothesis of the theorem are
  satisfied by (\ref{eq:ord_AQ}). Therefore, every Willmore sphere in
  $S^4$ either comes from a holomorphic curve in $\CP^3$ or is an
  inverted minimal sphere in $\R^4$ which, for immersed $f$, recovers
  the results by \cite{Bryant}, \cite{Ejiri}, \cite{montiel}.
\end{rem}

\begin{rem}
  As already mentioned in the introduction, there is evidence that the
  theorem holds without assuming the existence of a dual Willmore surface.
    But the proof of this conjecture, even in the genus 1
  case, seems more involved then simple degree computations. This is
  mainly due to the fact that the B\"acklund transform generally does
  not admit a mean curvature sphere congruence.
\end{rem}

\begin{proof}
  We may assume that $A\not= 0 $ and $Q \not= 0$.  In this case $AQ=0$
  implies that the forward and backward B\"acklund transforms
  coincide, i.e. $\tilde L = \hat L$. Moreover, $-S $ is the mean
  curvature sphere congruence of $\hat{L}$ since  
\[
\hat Q|_{\hat{L}} = 
  A|_{\tilde L} =0\,.
\]
Therefore $\tilde{L}$ is a dual Willmore surface and
our aim is to show that $\tilde L$ is in fact a point on the Willmore
surface $f$. Since $S$ stabilizes $\tilde L$, all the mean curvature
spheres of $f$ will then pass through a common point.  Inverting $f$
at this point thus gives a minimal surface in $\R^4$.

Assuming that $\tilde L$ is not a point its derivative
$\tilde\delta\in H^0(K\Hom_+(\overline{\tilde L}, \overline{V/\tilde
  L}))$ is a non--trivial holomorphic bundle map
(\ref{eq:delta_holomorphic}), so that its vanishing order calculates to
\begin{equation}
\label{eq:ord_delta_+}
\ord\tilde\delta=\ord\hat\delta= \deg K - \deg V + 2 \deg \tilde L\,.
\end{equation}
From (\ref{eq:ord_delta}),  (\ref{eq:ord_A}) and (\ref{eq:ord_delta_+}),
we obtain
\[ 
0 \le 2\ord A + \ord\delta+\ord\tilde \delta = 4\deg K + 2\deg
V\,,
\]
and similarly (\ref{eq:ord_delta}), (\ref{eq:ord_Q})  and
(\ref{eq:ord_delta_+}) give
\[
0 \le 2\ord Q + \ord\delta+\ord\hat \delta = 4\deg K - 2\deg
V\,.\]
Therefore, 
\[
|\deg V |\le 2\deg K\,,
\]
which  contradicts the degree  assumption of the theorem.
\end{proof}

\section{Loops of flat connections}
In addition to degree estimates, we now study the monodromies of the
associated family of flat connections arising from a Willmore surface
$f: M \to S^4$. The main reference for this is Section 6 of
\cite{Klassiker}. Recall (\ref{eq:Willmore}) that $f: M \to S^4$ is
Willmore if and only if
\[ 
d^\nabla * A =0\qquad \text{ or, equivalently, }  \qquad d^\nabla*Q=0\,,
\]
where $\nabla S = 2(*Q -*A)$ is the derivative of the mean curvature
sphere congruence $S: M \to \Zz$ of $f$.  One can immediately verify
that these equations are equivalent to the flatness of the family of
quaternionic connections
\begin{equation}
\label{eq:lambda_family_A}
\nabla_\lambda = \nabla + (\lambda -1) A
\end{equation}
where $\lambda=\alpha + \beta S$ with $\alpha,\beta\in \R$ and
$\alpha^2 + \beta^2 = 1$. The geometric interpretation of this family
of connections is the following: viewing the line bundle $L\subset V$
corresponding to $f$ in the flat background connection
$\nabla_\lambda$, we obtain the associated family of Willmore surfaces
$f_\lambda$ which generally have M\"obius monodromy.

For our purposes it is advantageous to extend $\nabla_\lambda$ to a
holomorphic family of flat complex connections parameterized over
$\C_*$.  To do this, we view $V$ as a complex vector bundle with
respect to the complex structure $I$ given by multiplication $I \psi =
\psi i$ by the quaternion $i$.  Then
\[
\lambda = \frac{\mu + \mu\invers}{2} + \frac{\mu\invers - \mu }{2}IS\,,
\]
where $\mu = a + I b\in\C_*$, extends $\lambda$ away from the unit
circle and  $\nabla_\lambda$ becomes
\begin{equation}
\label{eq:mu_family}
 \nabla_\mu =  (\nabla - A) + (\plusproj \mu + \minusproj
\mu\invers) A\,.
\end{equation}
Since $I$ is parallel with respect to $\nabla$, we see that the
flatness of the family $\nabla_\lambda$ is equivalent to the flatness
of the holomorphic family of complex connections $\nabla_\mu$ for
$\mu\in\C_*$. 

It is important to notice that the  $(0,1)$--part
with respect to the complex structure $S$ of 
$\nabla_\mu$ is independent of $\mu\in\C_*$ and gives the quaternionic holomorphic structure 
\begin{equation}
\label{eq:01_nabla_mu}
(\nabla_\mu)'' =\nabla''= \delbar + Q\,.
\end{equation}
In particular, every parallel section of
$\nabla_\mu$ for some $\mu\in\C_*$ is holomorphic, i.e., contained in $H^{0}(V,\delbar + Q)$.
We denote by
\begin{equation}
\label{eq:monodromy_mu}
H_\mu: \pi_1(M)\to \Gl(4,\C), \quad \mu\in \C_*,
\end{equation}
the holomorphic family of monodromy representations of the flat
connections $\nabla_\mu$. Notice that for unitary $\mu$ the connection
$\nabla_\mu$ is quaternionic and therefore
\begin{equation}
\label{eq:monodromy_lambda}
H_\mu: \pi_1(M)\to \Gl(2,\H), \quad \mu\in S^1\,.
\end{equation}

For a Willmore torus $f: T^2\to S^4$ with non--trivial normal bundle
the monodromies of the holomorphic family of complex connections
$\nabla_\mu$, and thus also the monodromies of the associated family
of Willmore surfaces, are either all trivial or translational.  
\begin{lemma}
\label{l:eigenvalue_one}
Let $f: T^2 \to S^4$ be a Willmore torus with non--trivial normal
bundle where $T^2 = \R^2/\Gamma$. Then $1$ is the only occurring
eigenvalue for the holomorphic family of monodromy representations
$H_\mu: \Gamma\to\Gl(4,\C)$.
\end{lemma}
\begin{proof}
Note that 
\[
\nabla_\mu =   (\nabla - A) + (\plusproj \mu + \minusproj
\mu\invers) A
\]
for $\mu = e^{tI}, \ t\in\C$, is gauge equivalent by $e^{\frac{t}{2}S}$ 
to 
\[
\tilde\nabla_\mu =   (\nabla - Q) + (\plusproj \mu\invers  + \minusproj
\mu) Q\,.
\]
Therefore, if $A=0$ or $Q=0$ the monodromy representation $H_\mu$ is
trivial for all $\mu\in\C_*$.

We now assume that $A\not=0$, $Q\not=0$, and that there exists
$\gamma\in\Gamma$ so that the family $H_\mu(\gamma)$ has non--constant
eigenvalues $h_\mu$ depending holomorphically  on $\mu$. In other
words, there exists a $\nabla_\mu$ parallel section
$\psi_\mu\in\Gamma(\pr^*V)$ on the universal cover $\pr:\R^2\to T^2$
with $H_\mu(\gamma)\psi_\mu = \psi_\mu h_\mu$. From
(\ref{eq:01_nabla_mu})  we see that the quaternionic holomorphic
structure on $V/L$ satisfies
\[
\delbar + \pi Q = \pi \nabla'' = \pi \nabla_\mu''\,. 
\]
In particular, $\varphi_\mu = \pi\psi_\mu \in \Gamma(\pr^*(V/L))$ is a
quaternionic holomorphic section with mono\-dromy $h_\mu$, i.e.,
\[ 
(\delbar + \pi Q)\varphi_\mu=0, \quad 
\gamma^*\varphi_\mu =\varphi_\mu h_\mu\,.
\]
First, we note that $\varphi_\mu\not=0$ since otherwise $\psi_\mu$
would be a parallel section of $\pr^*L$ which would imply that $f$ is
constant. Second, since $h_\mu$ is a non--constant holomorphic
function of $\mu$, the sections $\varphi_\mu\in\Gamma(\pr^*(V/L))$ are
linearly independent for $\mu$ near $\mu_0$ with $h'_{\mu_0}\not=0\,$:
the $\varphi_\mu$ are eigenvectors with distinct eigenvalues $h_\mu$
of the deck transformation operator $\gamma^*: \Gamma(\pr^*(V/L)) \to
\Gamma(\pr^*(V/L))$. On the other hand, the Pl\"ucker formula
\cite{Klassiker} for holomorphic sections with monodromy of the
quaternionic holomorphic line bundle $V/L$ bounds the number $n$ of
such linearly independent sections by its Willmore energy
\[
\mathcal{W}(V/L) \ge n((n-1)(1-g) -\deg(V/L)) = - n \deg (V/L)\,.
\]
Here $g$ is the genus of the underlying Riemann surface, which in our
case is $g=1$.

In case the normal bundle degree of $f$ satisfies $v =\deg V <0$, we
see from (\ref{eq:degV}) that also $\deg(V/L) < 0$.  Therefore, the
eigenvalues of $H_\mu(\gamma)$ must be independent of $\mu$ for all
$\gamma\in\Gamma$. Since $H_1$ is the trivial representation all the
eigenvalues of $H_\mu(\gamma)$ are equal to $1$. 

If the normal bundle degree of $f$ is positive, we will apply the
above argument to the Willmore surface $f^\perp: T^2 \to S^4$ given by the
line bundle $L^\perp\subset V^*$ whose mean curvature sphere
congruence is $S^*$ : since
\[
\nabla^*S^* = (\nabla S)^* = 2(-*\!A^* + *Q^*)\,,
\]
where $A^*\in\Gamma(\bar K\End_-(V^*))$ and
$Q^*\in\Gamma(K\End_-(V^*))$, we see that 
\[
Q^\perp = -A^*\,, \quad  A^\perp = -Q^*
\] 
and hence $L^\perp\subseteq \ker Q^\perp$.  Moreover,
$d^{\nabla^*}*Q^\perp =0$ so that $f^\perp$ is also Willmore.  The
corresponding family of flat connections is given by
\[
\nabla^\perp_\mu=  (\nabla^*-A^\perp) +
(\plusproj\mu + \minusproj \mu\invers)A^\perp\,
\] 
which, as we have seen above, is gauge equivalent to
\[
\tilde\nabla^\perp_\mu=  (\nabla^*-Q^\perp) +
(\plusproj\mu\invers + \minusproj \mu)Q^\perp\,.
\]
But the latter is the dual connection of $\nabla_\mu$ so that
$\nabla^\perp_\mu$ is gauge equivalent to $(\nabla_\mu)^*$.
Therefore, the monodromy representations $H_\mu$ and $H^\perp_\mu$
have the same eigenvalues.  If the normal bundle degree of $f$ is
positive, i.e., $v = \deg V >0$, then $V^*$ with complex structure
$S^*$ has negative degree $v^* = \deg V^* <0$ and we can apply our
previous argument to $f^\perp$. Again we deduce that all the
eigenvalues of $H_\mu^\perp$, and thus also of $H_\mu$, are equal to
$1$.
 \end{proof}

\begin{rem}
  In the previous proof, we used the Pl\"ucker formula for holomorphic
  sections with monodromy whereas in \cite{Klassiker} this formula is
  only proven for holomorphic sections without monodromy.  To allow
  for monodromy, we adapt the proof in \cite{Klassiker} to our
  situation by replacing the trivial connection with a flat connection.
\end{rem}

From the previous lemma,  we see that $V$ admits a
$\nabla_\mu$--parallel complex line subbundle $U_\mu\subset V$. For
$|\mu|=1$ the connection $\nabla_\mu$ is quaternionic and thus we
obtain a $\nabla_\mu$--parallel quaternionic line subbundle.

\begin{lemma}
\label{l:translations_holonomy}
  Let $V$ be a rank 2 quaternionic vector bundle over a torus
  $T^2=\R^2/\Gamma$ with flat connection $\nabla$.  Assume that the
  monodromy representation $H: \Gamma\to \Gamma(\Gl(V))$ of $\nabla$
  has  $1$ as its only  eigenvalue.
  
  Then there exists a parallel quaternionic line subbundle $U\subset
  V$ on which $\nabla$ is trivial.  If we denote by $\Rr:=\Hom(V/U,U)$,
  then
\[
R:=H-\Id: \Gamma\to\Gamma(\Rr)
\]
is a translational representation.  Moreover, there exists
$\omega\in\Omega^1(\Rr)$ of the form $\omega= B_1 dx + B_2 dy$ with $
B_i\in\Gamma(\Rr)$ parallel with respect to $\nabla$, such that $\nabla +
\omega$ is a trivial connection.

\end{lemma}


\begin{proof}
  Since $1$ is an eigenvalue of $H$ there exists $\psi\in\Gamma(V)$
  with $H\psi=\psi$. But $\psi$ is nowhere vanishing and thus spans a
  parallel quaternionic line subbundle $U$ on which $\nabla$ is
  trivial, i.e., $R|_U=0$.  On the other hand, the
  characteristic polynomial of $R$ is $X^4$ and hence $\tr_\C R^n =0$.
  This implies $R^2=0$ and therefore $RV\subset U$.  It is easy to
  check that $R_{\gamma_1\gamma_2} = R_{\gamma_1} + R_{\gamma_2}$ so
  that $R = H-I$ gives a representation into $\Gamma(\Rr)$. In
  particular, this implies that the induced connection $\nabla$ on
  $\Rr$ is trivial and that $\nabla R_\gamma=0$.  For fixed $p\in T^2$
  the representation $R(p): \Gamma\to\Rr_p$ is given by $R_\gamma(p) =
  \int_\gamma \omega(p)$, where $\omega(p)=B_1(p)dx+B_2(p)dy$ is an
  $\Rr_p$--valued harmonic form. Because $\nabla R_\gamma=0$, the
  sections $B_i\in\Gamma(\Rr)$ satisfy $\nabla B_i =0$ and hence the
  $\Rr$--valued 1--form $\omega\in\Omega^1(\Rr)$ is closed, i.e., 
  $d^\nabla\omega =0$. This implies that the connection
  $\nabla+\omega$ is flat. To see that $\nabla+\omega$ has no
  monodromy, we let $\varphi\in\Gamma(\pr^*V)$ be a $\nabla$--parallel
  section and define
 \[
 \tilde\varphi:=\varphi -(\int_{p_0}\omega)\varphi\,,
  \]
  where $p_0\in T^2$ is a chosen base point.  Then it is easy to check
  that $\tilde\varphi$ is parallel with respect to $\nabla+\omega$ and
  has no monodromy.
\end{proof}

\begin{cor}
\label{cor:hol_parallel}
In the situation of the previous lemma, we denote by 
\[
\delbari = \delbar_0 - \omega^{(0,1)}
\]
the holomorphic structure with respect to the complex structure $I$ on
$V$. Here $\delbar_0=\tilde\nabla^{(0,1)}$ denotes the trivial
holomorphic structure on $V$.  Then the holomorphic sections of
$\delbari$ are the parallel sections of $\nabla$ which, if
$\omega\not=0$, are contained in $U$. In particular, we have a 4 or
2--dimensional space of holomorphic sections depending on whether
$\omega=0$ or not.
\end{cor}

\begin{proof}
  Let $U_1=U\subset V$ be the quaternionic line subbundle on which
  $\nabla$ is trivial. Since $U\subset \ker\omega$ and $\tilde\nabla
  =\nabla + \omega$, we see that $U_1\subset V$ is also
  $\tilde\nabla$--trivial.  Let $U_2\subset V$ be a complementary
  $\tilde\nabla$--trivial subbundle so that $V = U_1\oplus U_2$. 
If $\varphi=\varphi_1+\varphi_2\in\Gamma(V)$ is a holomorphic
  section, i.e.,
\[
\delbari\varphi=  \delbar_0\varphi - \omega^{(0,1)} \varphi = 0\,,
\]
then the latter is equivalent to
\[
\delbar_0\varphi_1= \omega^{(0,1)} \varphi_2\,, \quad 
\delbar_0\varphi_2 = 0\,.
\]
From our assumptions, we see that $\omega^{(0,1)}=Bd\bar z$ where $B =
\tfrac{1}{2}(B_1 - I B_2)$ is parallel with respect to $\nabla$, and
hence also with respect to $\tilde\nabla$. Therefore, $\varphi_2$ is
$\tilde\nabla$--parallel and $\delbar_0\varphi_1= B\varphi_2 d\bar z$.
This implies that $\varphi_1$ is harmonic on the torus $T^2$ and thus
$\tilde\nabla$--parallel. If $\omega\not=0$ then $B\varphi_2=0$ shows
that $\varphi_2=0$. Since $\tilde\nabla = \nabla +\omega$, we see that
$\varphi$ is $\nabla$--parallel.

\end{proof}

\begin{rem}
\label{rem:dim_parallel}
The flat connections $\nabla_\mu$ for $|\mu|=1$ are quaternionic.
Thus, applying the previous corollary to the flat connections
$\nabla_\mu$ for $|\mu|=1$, we see that the number of holomorphic
sections of $\nabla_\mu^{(0,1)}$ is either 4 or 2, depending on
whether $\omega_\mu=0$ or not. Since the dependence on $\mu$ is
holomorphic, this holds also for $\mu\in\C_*$.  Therefore every
holomorphic section of $\nabla_\mu^{(0,1)}$ is parallel with respect
to $\nabla_\mu$ for $\mu\in\C_*$.
\end{rem}

\section{Willmore tori with  non--trivial normal bundle}

In the previous section we have seen that the monodromy of the
associated family of a Willmore torus $f: T^2 \to S^4$ with
non--trivial normal bundle is either trivial or translational. The
former occurs for the twistor projection of a holomorphic curve in
$\CP^3$ since in this case $A=0$ or $Q=0$. On the other hand,
translational monodromy occurs from the periods around the ends of the
associated family of a minimal surface in $\R^4$ with planar ends.
The main result of this paper is that these are in fact the only
possibilities:

\begin{theorem}
\label{th:main}
Let $f: T^2\to S^4$ be a Willmore torus with non--trivial normal
bundle.  Then, either $f$ or $f^\perp$ is a twistor projection of an
elliptic curve in $\CP^3$, or $f$ is an inverted minimal torus in
$\R^4$.
\end{theorem}
\begin{proof} We may assume that $A$ and $Q$ are not identically zero.
  Due to Theorem \ref{t:AQ=0} it suffices to show that under our
  assumptions $AQ =0$, i.e., that $f$ admits a dual Willmore surface. Assume this were not the case. Since we are
  working over a torus, (\ref{eq:ord_AQ}) shows that $AQ$ and $\delta$
  have no zeros.  In particular, $A$ and $Q$ have no zeros and
  therefore (\ref{eq:ord_A}), (\ref{eq:ord_Q}) imply that
\begin{equation*}
\deg \tilde L =  3\deg L, \quad \deg \hat L =-\deg L\,.
\end{equation*}
We may assume that the normal bundle degree $v$ of $f$ is positive.
Otherwise we work with the Willmore surface $f^\perp$. Therefore,
(\ref{eq:degV}) and (\ref{eq:normalbundle_degree}) imply that $\deg L
=\tfrac{1}{2} v >0$ and hence
\begin{equation}\label{eq:deg_kerA}
\deg\tilde L > 0\,.
\end{equation}

Taking the $(1,0)$--parts of the complex connections
\begin{equation}
 \nabla_\mu =  (\nabla - A) + (\plusproj \mu + \minusproj
\mu\invers) A
\end{equation}
with respect to the complex structure $I$, gives the holomorphic family
of antiholomorphic structures
\[
\partial_\mu:= \nabla_\mu^{(1,0)} = \nabla^{(1,0)} + (\mu-1)\plusproj A
\]
on the complex vector bundle $V$. Here we have
used that $A^{(1,0)} = \plusproj A$.

Even though the holomorphic family of flat connections $\nabla_\mu$
does not extend into $\mu=0$, the family of antiholomorphic structures
$\partial_\mu$ does.  Corollary \ref{cor:hol_parallel} and Remark
\ref{rem:dim_parallel} show that every antiholomorphic section is
$\nabla_\mu$--parallel for $\mu\not=0$ and thus by (\ref{eq:01_nabla_mu})
holomorphic with respect to $\nabla_\mu'' = \delbar+Q$. In particular, the kernels of $\partial_{\mu}$ are all contained in the finite dimensional vector space $H^0(V,\delbar+Q)$ of quaternionic holomorphic sections of $V$. Consider
\[
\partial_\mu: H^0(V,\delbar+Q) \to \Omega^{(1,0)}(V)
\]
as a holomorphic family of endomorphisms with finite dimensional
domain parameterized over $\mu\in\C$. Then the minimal kernel dimension
of $\partial_\mu$ is generic, and we obtain a complex holomorphic
vector bundle $\K$ of rank 2 or 4 over $\C$ with
$\K_\mu\subseteq\ker\partial_\mu$.  If $\psi$ is a holomorphic
section of $\K$, then
\[
\psi(\mu)= \psi_0+ \mu\psi_1 + O(\mu^2)
\]
is parallel with respect to $\nabla_\mu$ for $\mu\not=0$ and
$\psi_0\in\K_0\subset H^0(V, \delbar+Q)$.

Recall the decomposition  $V = V_+ \oplus V_-$ into the $\pm i$
eigenspaces $V_\pm = \tfrac{1}{2}(1 \mp IS)V$ of $S$.
Since $(1\pm IS)A = A(1\mp IS)$, we obtain
\[
\nabla^\mu\psi_\mu =  \mu\invers A\psi_0^+ + (\nabla -A)\psi_0 + A\psi_1^+ +
 O(\mu) =0
\]
for $\mu\in\C_*$. Comparing coefficients at powers of $\mu$ gives
\[
A\psi_0^+ = 0 \quad \text{ and } \quad (\nabla - A)\psi_0 + A\psi_1^+ = 0\,.
\]
Finally, taking $(1,0)$ and $(0,1)$--parts with respect to the complex
structure $S$, we arrive at
\begin{equation}
\label{eq:final}
 A\psi_0^+ = 0\,,  \quad \partial \psi_0^+ = 0\,, \quad  \partial
 \psi_0^- + A\psi_1^+ = 0\,, \quad \text{ and } \quad (\dbar + Q) \psi_0 = 0\,,
\end{equation}
where we again used the direct sum decomposition $V = V_+ \oplus V_-$.
The first two equations of (\ref{eq:final}) imply that $\psi_0^+$ is
an antiholomorphic section of the complex line bundle $\tilde L_+$, where $\tilde L=\ker A$ is
the forward B\"acklund transform (\ref{eq:baecklund}) of $f$. But
$\tilde L$ has positive degree by (\ref{eq:deg_kerA}) and therefore
$\psi_0^+=0$. Decomposition of the last equation in (\ref{eq:final})
according to $V=V_+ \oplus V_-$ gives
\[
\dbar\psi_0^-=0 \quad \text{ and } \quad Q\psi_0^- = 0\,.
\]
We now recall that $L =\ker Q$ so that $\psi_0^-\in\Gamma(L)$ is a
section of $L$ which, by (\ref{eq:final}), satisfies
\[
\nabla\psi_0^-= A(\psi_0^--\psi_1^+)\,.
\]  
Since $A$ has image in $L$ and $\psi_0^-$ can be chosen not 
identically zero, this implies that $L\subset V$ is
$\nabla$--parallel, i.e., that $f$ is constant.
\end{proof}
\bibliographystyle{amsplain}

\end{document}